\documentclass[a4paper]{article}

\usepackage{amsmath,amsfonts,latexsym,theorem
}
\nonstopmode

\numberwithin{equation}{section}

\newtheorem{Theorem}{Theorem}[section]

\newtheorem{Proposition}{Proposition}[section]
\newtheorem{Lemma}{Lemma}[section]

\newtheorem{Corollary}{Corollary}[section]
\newenvironment{Proofc}[1]{\smallskip\par\noindent\textsc{#1}\quad}%
  {\hfill$\Box$\bigskip\par}
\newenvironment{Proof}{\begin{Proofc}{Proof}}{\end{Proofc}}
\newenvironment{Remc}[1]{\smallskip\par\noindent\textsc{#1}\quad}{\smallskip\par}
\newtheorem{Remark}{Remark}[section]

\def\a{\alpha}

\def\d{\delta}
\def\D{\Delta}
\def\g{\gamma}
\def\l{\lambda}

\def\s{\sigma}
\def\o{\omega}
\def\th{\theta}
\def\Th{\Theta}
\def\t{\tau}
\def\e{\varepsilon}
\def\z{\zeta}

\newcommand{\hx}{{\hat x}}

\newcommand{\Ho}{{\overline H}}
\newcommand{\xo}{{\overline x}}
\newcommand{\tov}{{\overline t}}
\newcommand{\po}{{\overline p}}
\newcommand{\uo}{{\overline u}}
\newcommand{\Xo}{{\overline X}}

\newcommand{\tr}{\operatorname{\text{tr}}}
 \newcommand{\Z}{{\mathbb Z}}
\newcommand{\R}{{\mathbb R}}
\newcommand{\Rn}{{\mathbb R^n}}

\newcommand{\E}{{\mathbb E}}

\newcommand{\tz}{{\th\z}}

\begin{document}
\title{Continuous dependence estimates and homogenization of  quasi-monotone systems of fully nonlinear second order parabolic  equations}

\author{Fabio Camilli\footnotemark[1] \and Claudio Marchi\footnotemark[2]}

\date{version: \today}
\maketitle

\footnotetext[1]{Dip. di Scienze di Base e Applicate per l'Ingegneria,  ``Sapienza" Universit{\`a}  di Roma, via Scarpa 16,
 00161 Roma, Italy, ({\tt e-mail:camilli@dmmm.uniroma1.it})}
\footnotetext[2]{Dip. di Matematica, Universit\`a di Padova, via Trieste 63, 35121 Padova, Italy ({\tt marchi@math.unipd.it}).}

\begin{abstract}
Aim of this paper is to extend  the  continuous dependence estimates proved in  \cite{JK1} to quasi-monotone systems
of  fully nonlinear second order parabolic equations. As by-product of these estimates, we get   an H\"older   estimate for
bounded  solutions of  systems and a rate of convergence  estimate for the  vanishing viscosity approximation.

In the second part of the paper we employ similar techniques
to study the periodic homogenization of quasi-monotone systems of fully nonlinear second order uniformly parabolic equations. Finally, some examples are discussed.
\end{abstract}

 \begin{description}
\item [{\bf MSC 2000}:] 35J88, 49L25, 35B27.
 \item [{\bf Keywords}:] Continuous dependence estimate, quasi-monotone system,   Hamilton-Jacobi-Bellman-Isaacs equation, vanishing viscosity, ergodic equation, homogenization.
 \end{description}


\section{Introduction}\label{intro}
This paper is devoted to the weakly coupled  system of parabolic equations
\begin{equation}\label{P}
\partial_t u_i+H_i(t,x,u,Du_i,D^2u_i)=0 \qquad\textrm{in }(0,T)\times\R^n,\, i=1,\dots,m
\end{equation}
where $\partial_t\equiv \partial/\partial t$, the operators $H_i:(0,T)\times\R^n\times\R^m\times\R^n\times S^n$ are given by
\begin{equation}\label{minmax}
H_i(t,x,r,p,X)=\min\limits_{\z\in Z_i}\max\limits_{\th\in \Th_i}
\left\{
-\tr\left(A_i^\tz(t,x,p)X\right)+f^\tz_i(t,x,r,p,X)\right\}
\end{equation}
and  $u(x)=(u_1(x),\dots, u_m(x))$.
In fact, our techniques may be easily adapted to the case of systems of  elliptic equations. Here, all (sub-, super-) solutions will always be in {\it viscosity} sense (see below for the precise definition; for the main properties, we refer the reader to \cite{IK91} and also to \cite{CIL} for a single equation).\par
Quasi-monotonicity is a basic assumption  which guarantees the validity
of the maximum principle   for   weakly coupled systems. In \cite{IK91}
this assumption has been exploited to prove general existence and uniqueness
results for solutions of  systems of fully nonlinear second order PDEs.
Aim of this paper is to show that this  assumption allows to extend to weakly coupled  systems two well known  properties of fully nonlinear equations:
continuous dependence estimates and periodic homogenization.\par
Continuous dependence estimates (namely, an estimate of $|u(t,x)-v(t,x)|$ where $u$ and $v$ are two solutions to \eqref{P}-\eqref{minmax} with different coefficients) are  useful tools to obtain regularity results and rate of convergence estimates (e.g. for vanishing viscosity and numerical approximation). A general result  for, possibly degenerate, scalar equations
was proved in \cite{JK1,JK2} (see also \cite{m}) using techniques based on the maximum principle for semi-continuous solutions: doubling the variables and adding a penalization term.
We show that the quasi-monotonicity assumption allows to extend the result in \cite{JK1} to weakly coupled systems at the same level of generality
(see also \cite{BJK}, \cite{bcz} for related results).
As an application of continuous dependence estimates, we obtain regularity estimates (a priori $L^\infty$ and  H\"older bounds; we refer the reader to \cite{BS, FP} for Harnack type estimates for these systems) and a rate of convergence estimate for the vanishing viscosity approximation (in this direction, this paper extends the results in~\cite {BCD} to the case of quasi-monotone systems).
We shall also illustrate our results for a class of systems which arises in optimal control theory and, as scalar equation, it encompasses the Hamilton-Jacobi-Bellman-Isaacs equation associated to stochastic differential games (see \cite{FSo, FSu}). In this case, taking advantage of the special form of the coefficients, we obtain a simpler expression of the estimates.

\vskip 1mm

In the second part of the paper we are concerned with periodic homogenization  of weakly coupled systems of uniformly parabolic equations. In this case, the coefficients in \eqref{minmax} have the form
\begin{equation*}
A_i^\tz=A_i^\tz\left(x,\frac x\e\right), \qquad f^\tz_i=f^\tz_i\left(x,\frac x \e,r,p\right)
\end{equation*}
and are $\Z^n$-periodic in the $x/\e$-variables. The parameter $\e$ is meant to tend to $0$. This fact modelizes a medium displaying heterogeneities in a microscopic scale while one seeks a description of the macroscopic phenomena (which are the only relevant ones). At the limit, the solutions are expected to converge to the solution of a ``homogenized'' problem where the effective operator needs to be suitably defined.

The homogenization problem for a scalar equation has been studied e.g. in \cite{AB2,AL,Ev2,hi,LPV} (see \cite{AB8} for a general review of the results). A homogenization result
for  quasi-monotone systems of first order Hamilton-Jacobi equations was obtained in \cite{cll}. For systems of second order equations we refer the reader to \cite{BLP, B05} for the quasi-linear case  and to \cite{BH,BHP} for homogenization via probabilistic techniques. Also for the homogenization we are able to generalize the result from   the scalar case to the weakly coupled one making a crucial use of the quasi-monotonicity assumption.
The proof   relies on an appropriate modifications of the perturbed test function's method introduced by Evans \cite{Ev2}.

This paper is organized as follows: in the rest of this section, we introduce our notations, we list the standing assumptions and we recall the definition of viscosity solution of a system of PDEs.
Section \ref{results} is devoted to the continuous dependence estimate; in particular, we illustrate our results for a class of systems which arise in optimal control problems.
Taking advantage of this estimate, in Section \ref{reg_vv}, we deduce a regularity result and a rate of convergence for the vanishing viscosity; moreover, we work out in detail the vanishing viscosity approximation of a first order system arising in optimal control.
Section \ref{homogenization} is devoted to homogenization results. Finally, in the Appendix we give the proof of a technical Lemma and, for the sake of completeness, we quote some results yet established in the literature.


\subsection{Notations and standing assumptions}\label{assumption}

\noindent \textbf{Notations: }
We set $I:=\{1,\dots, m\}$ and $Q_t:=(0,t)\times\R^n$. $S^n$ denotes the set of $n\times n$ real symmetric matrices; it is endowed with the Frobenius norm and the usual order, namely: $|X|=\tr(XX^T)^{1/2}$ and $X\geq Y$ whenever $X-Y$ is a semidefinite positive matrix. For each function~$h$ defined on $(0,T)\times \R^n$, $\bar{{\cal P}}^{2,+}h(\tau,\xi)$ and $\bar{{\cal P}}^{2,-}h(\tau,\xi)$ denote respectively the parabolic super- and subjects at the point $(\tau,\xi)$ (see \cite[Section~8]{CIL}). For $f:\R^n\to\R^m$, we define the $C^0$-norm by $\|f\|=\sup_{i\in I,\,x\in\R^n} |f(x)|$  and, for $\mu\in(0,1]$, the H\"older seminorm by $[f]_\mu:=\sup_{i\in I,\,x\ne y}\frac{|f_i(x)-f_i(y)|}{|x-y|^\mu}$.
For $\mu\in(0,1]$, $C^\mu(\R^n)$ denotes the H\"older space of functions~$f$ such that: $\|f\|+[f]_\mu <+\infty$. Finally, $BUC(\R^n)$ denotes
 the space of uniformly continuous, bounded functions $f:\R^n\to\R^m$.

\vskip 1mm

\noindent\textbf{Standing assumptions:}
For $i\in I$ and $H_i$ defined as in \eqref{minmax}, we assume
\begin{itemize}
\item[($C0$)] The sets $\Th_i$ and $Z_i$ are compact metric spaces. Moreover, wlog, we assume $\Th_i=\Th$ and $Z_i=Z$ (it suffices to consider $\Th=\Pi_i \Th_i$, $Z=\Pi_i Z_i$ and to  extend the functions $f_i$ and $A_i$ to these sets).
\item[($C1$)]   For every $R>0$, $f^{\tz}_i\in C([0,T]\times\R^n\times\R^m\times\R^n\times S^n)$ is uniformly continuous on the set $[0,T]\times\R^n\times[-R,R]^m\times B_{\R^n}(0,R)\times B_{S^n}(0,R)$ uniformly in $\th,\z$.
\item[($C2$)] For every $X,Y\in S^n$ with $X\le Y$ there holds
\[f^\tz_i(t,x,r,p,X)\ge f^\tz_i(t,x,r,p,Y)\qquad \forall t,x,r,p,\th,\z.\]
\item[($C3$)]  For every $t,x,p,\th,\z$, $A^\tz_i(t,x,p)=a^\tz_i(t,x,p)a^\tz_i(t,x,p)^T$ for some matrix $a^\tz_i\in C([0,T]\times \R^n\times
\R^n).$ Furthermore, for every $R>0$, $a^\tz_i$ is uniformly continuous on $[0,T]\times\R^n\times B_{\R^n}(0,R)$ uniformly in $\th,\z$.
\item[($C4$)]
For every $R>0$, there is $\g_R\in\R$ s.t. if  $r$, $s\in[-R,R]^m$
and $\displaystyle r_j-s_j=\max_{k\in I}\{r_k-s_k\}\ge 0$, then
\[f^\tz_j(t,x,r,p,X)-f^\tz_j(t,x,s,p,X)\ge\gamma_R(r_j-s_j) \qquad \forall t,x,p,X,\th,\z. \]
\item[$(C5)$] There exists $\mu\in(0,1]$ such that: for every $R>0$ there exists a constant $C_{f,R}$ such that
$$
\left| f^\tz_i(t,x,r,p,X)-f^\tz_i(t,y,r,p,X)\right|\leq C_{f,R}\left(|p||x-y|+|x-y|^\mu\right)
$$
for every $\th,\z,i,t,x,y,p,r,X$ with $|r|<R$.
\item[$(C6)$]
There is a constant $C_a$ such that
$$
\left| a^\tz_i(t,x,p)-a^\tz_i(t,y,p)\right|\leq C_a|x-y| \qquad
\forall \th,\z,i,t,x,y,p.
$$
\item[$(C7)$] There holds: $C^f:=\sup_{\th,\z,i,t,x}|f^\tz_i(t,x,0,0,0)|\leq +\infty$.
\end{itemize}
\begin{Remark}
Assumption ($C4$)  implies a quasi-monotonicity property of the system~\eqref{P}; namely, for every $R>0$, there is $\g_R\in\R$ s.t. if  $r$, $s\in[-R,R]^m$
and $\displaystyle r_j-s_j=\max_{k\in I}\{r_k-s_k\}\ge 0$, then
 \begin{equation}\label{QM}
  H_j(t,x,r,p,X)-H_j(t,x,s,p,X)\ge\gamma_R(r_j-s_j) \qquad \forall t,x,p,X.
 \end{equation}
\end{Remark}
\begin{Remark}
We refer the reader to Section~\ref{oc_pb} for a class of systems (arising in optimal control theory) which fulfills assumptions $(C0)$-$(C4)$.
Let us also observe that, when system \eqref{P}-\eqref{minmax} reduces to a single equation, the above assumptions are satisfied, e.g., by: the Hamilton-Jacobi-Bellman-Isaacs equation associated to a two-players zero-sum stochastic differential game, the equation of mean curvature flow of graphs, the $p$-Laplacian with $p>2$ (see \cite{CIL, JK1}).
Furthermore, let us recall that a wide class of nonlinear operators can be written in the form \eqref{minmax} (see \cite{Ka95, ES84}).\par
\end{Remark}

\vskip 1mm

\noindent\textbf{Definition of solution (\cite{IK91}):}
\noindent($i$) An USC function $u:Q_T\rightarrow \R^m$ is a {\it subsolution} of \eqref{P} if: whenever $\phi\in C^2(Q_T)$, $i\in I$ and $u_i-\phi$ attains a local maximum at $(t,x)$, then there holds
\begin{equation*}
\partial_t \phi(t,x) +H_i(t,x,u(t,x),D\phi(t,x),D^2\phi(t,x))\leq 0. 
\end{equation*}

\noindent($ii$) A LSC function $u:Q_T\rightarrow \R^m$ is a {\it supersolution} of \eqref{P} if: whenever $\phi\in C^2(Q_T)$, $i\in I$ and $u_i-\phi$ attains a local minimum at $(t,x)$, then there holds
\begin{equation*}
\partial_t \phi(t,x) +H_i(t,x,u(t,x),D\phi(t,x),D^2\phi(t,x))\geq 0. 
\end{equation*}

\noindent($iii$) A function $u$ is a {\it solution} of \eqref{P} if it is both a sub- and a supersolution. In particular, it belongs to $C(Q_T)$.



\section{The continuous dependence estimate}\label{results}
In this section we prove the continuous dependence estimate for the problem \eqref{P}-\eqref{minmax}.
%
%

\begin{Theorem}\label{thm:dipcnt}
Assume that, for $k=1,2$, $H^k=\{H_i^k\}_{i\in I}$ satisfies assumptions ($C0$)-($C4$) with constant $\g_R^k$.
Let $u^1$ and $u^2$ be respectively a bounded subsolution to problem~\eqref{P}-\eqref{minmax} with $H=H^1$
and a bounded supersolution to problem~\eqref{P}-\eqref{minmax} with $H=H^2$. Set  $R:=\max (\|u^1\|,\|u^2\|)$
and $\g=\min(\g_R^1,\g_R^2)$. Then for each $0\le t\le T$, $\bar\g\ge 0$ and $\a>0$, we have
\begin{align*}
&  \sup_{E_t^\a}\Big(e^{\g \t}\big(u^1_i(\t,x)-u^2_i(\t,y)\big) -\frac \a 2e^{\bar  \g \t}|x-y|^2\big)\Big)\leq\\
&\quad \sup_{E_0^\a}  \big(u^1_i(0,x)-u^2_i(0,y) -\frac \a 2|x-y|^2\big)^+  + t\sup_{D^\a_{\g t}} \Big(e^{\g \t}[f^{\tz,1}_{i}(\t,y,r,p,X)-f^{\tz,2}_{i}(\t,x,r,p,X)]\\
&\qquad +3\a e^{\bar \g\t}|a^{\tz,1}_{i}(\t,x,p)-
a^{\tz,2}_{i}(\t,y,p)|^2-\frac \a 2\bar\g e^{\bar \g \t}|x-y|^2\Big)^+
\end{align*}
where
\begin{align}
    &\Delta^\a:=\Big\{(x,y)\in\R^n\times \R^n:\,|x-y|\le 2\frac{R^{1/2}}{\sqrt \a} \Big\}\label{Deltaalpha}\\
    &E^\a_t:=\{(\t,x,y,i): 0\le \t\le t,\, (x,y)\in\Delta^\a, \, i\in I\}\nonumber\\
     \begin{split}\nonumber
     &D^\a_{\g t}:=\{(\t,x,y,i,r,p,X,\th,\z):\, p=\a(x-y)e^{(\bar\g-\g)\t},\,
     (\t,x,y,i)\in E^\a_t,\\
     & \qquad\qquad\quad|r|\le e^{-\g t}\min(\|u^1\|,\|u^2\|), \,|X|\le 3\a n e^{(\bar\g-\g)\t}, \th \in \Th, \z\in Z\}.
     \end{split}
\end{align}
\end{Theorem}
\begin{Proof}
We first consider the case $\g=0$.
Without loss of generality, assume $\|u_1\|\leq \|u_2\|$ (the other case can be dealt with in a similar manner and we shall omit it). Fix $t\in(0,T]$, $\a>0$ and $\bar \g\geq 0$. For every $0<\e\leq \a/5$, we set
\begin{eqnarray*}
\s^0&:=&\sup_{E^\a_0}\big(u^1_i(0,x)-u^2_i(0,y)-\frac\a 2|x-y|^2\big)^+\\
\s&:=&-\s_0+\sup_{E^\a_t}\big\{u^1_i(\t,x)-u^2_i(\t,y)-\big(\frac \a 2 e^{\bar \g \t} |x-y|^2
+\frac \e2(|x|^2+|y|^2)+\frac {\e}{t-\tau}\big)\big\}.
\end{eqnarray*}
Since we want to derive  an upper bound of $\s$, it is not restrictive to assume $\s>0$. For $\d\in (0,1)$, set
\begin{equation}\label{T4}
\psi(\t,x,y,i):=u_i^1(\t,x)-u_i^2(\t,y)-\frac{\d \s \t}{t}-
\left(\frac \a 2 e^{\bar \g \t}|x-y|^2+\frac \e 2(|x|^2+|y|^2)
+\frac {\e}{t-\tau}\right)
\end{equation}
for every $\t\in(0,t)$, $x,y\in\R^n$ and $i\in I$.
Since the functions $u_i^1$ and $u_i^2$ are bounded in~$Q_t$ and $\psi$ tends to $-\infty$ both as $\t\to t^-$ and as $|x|+|y|\to+\infty$, we deduce that there exists a point~$(\t_0,x_0,y_0, i_0)$ where the function~$\psi$ attains its global maximum, i.e.
\[
\psi(\t_0,x_0,y_0,i_0)\ge\psi(\t ,x ,y ,i ) \qquad \forall (\t,x,y,i)\in [0,t)\times\R^n\times\R^n\times I.
\]
By its definition \eqref{T4}, the function $\psi$ satisfies
\begin{equation}\label{T5}
\sup\limits_{E^\a_t}\psi
 \geq  \s+\s_0-\d\s=(1-\d)\s+\s_0.
\end{equation}

\begin{Lemma}\label{Lemma2}
Let $(\t_0,x_0,y_0,i_0)$ be the point where the function $\psi$ in \eqref{T4} attains its maximum. Then
\begin{itemize}
\item[$i)$] There holds
\begin{equation}\label{T7}
   |x_0-y_0|\leq 2\left(\frac{R}{\a}\right)^{1/2},\qquad  |x_0|, \,|y_0| \leq 2\left(\frac{R}{\e}\right)^{1/2}
\end{equation}
where $R$ is the constant introduced in Theorem \ref{thm:dipcnt}; in fact, there exists a modulus of continuity $m$ such that
\begin{equation}\label{T8}
   |x_0|,|y_0|\le \e^{-1/2}m(\e).
\end{equation}
\item[$ii)$] Assume that $u^1$ and $u^2$ are continuous in $x$ uniformly in $t$, namely, there exists a modulus of continuity $\o$ such that: $|u^j(\t,x)-u^j(\t,y)|\leq \o(|x-y|)$ ($j=1,2$). Then, we have
\begin{equation}\label{T7-1}
\a e^{\bar \g \t_0}|x_0-y_0|^2\leq \o(|x_0-y_0|).
\end{equation}
\item[$iii)$]
Assume that either $u^1$ or $u^2$ belongs to $C^1$. Then, we have
\begin{equation}\label{T7-2}
\a e^{\bar \g \t_0}|x_0-y_0|\leq n\left[\min_{j=1,2}\{[u^j]_1\} +\e^{1/2}\sqrt{2R}\right].
\end{equation}
\end{itemize}
\end{Lemma}
The proof is postponed to the Appendix.
We continue with the proof of Theorem \ref{thm:dipcnt}.

By Lemma~\ref{Lemma2}-(i), we deduce that $\t_0>0$; actually, for $\t_0=0$, inequality \eqref{T5} implies
\[\s_0+(1-\d)\s\le \psi(0,x_0,y_0, i_0)\le \s_0\]
and, in particular, $\s\le 0$, a contradiction.

We introduce the test function
\begin{equation*}
\phi(\t,x,y):=
\frac{\d \s\t}{t}+\frac \a 2 e^{\bar \g \t}|x-y|^2+\frac\e 2(|x|^2+|y|^2)+\frac {\e}{t-\tau}
\end{equation*}
and, for $i=i_0$ fixed, we invoke \cite[Thm 8.3]{CIL}: for every $\nu>0$, there exist values $a,b\in\R$ and matrices $X,Y\in S^n$ such that
\begin{align}\label{cil1}
\begin{split}
\left(a,p_{x_0},X\right)\in \bar{\cal P}^{2,+}u^1_{i_0}(\t_0,x_0),\qquad
\left(b,p_{y_0},Y\right)\in \bar{\cal P}^{2,-}u^2_{i_0}(\t_0,y_0),
\end{split}\\[8pt]
&
a-b=\partial_\tau\phi(\t_0,x_0,y_0)\equiv \frac{\d \s}t+\frac {\e}{(t-\tau_0)^2}+\frac{\a}{2}\bar\g e^{\bar \g\t_0}|x_0-y_0|^2\label{cil2} \\[8pt]
&-(\nu^{-1} +\bar \a+\e)\left(\begin{array}{cc} I&0\\0&I\end{array}\right)\leq\left(\begin{array}{cc} X&0\\0&-Y\end{array}\right)\leq
\Phi+\nu \Phi^2, &  \label{cil3}
\end{align}
where
\[\bar\a:=e^{\bar \g\t_0}\a,\quad \Phi:=\left(\begin{array}{cc}  \phi_{xx}&\phi_{xy}\\\phi_{yx}&\phi_{yy}\end{array}\right)_{(\t_0,x_0,y_0)}, \quad
p_{x_0}:=D_x\phi(\t_0,x_0,y_0),\quad p_{y_0}:=-D_y\phi(\t_0,x_0,y_0)
\]
(note that, according to notations of \cite{CIL}, the norm of a symmetrix matrix $A$ is defined as follows: $|A|_*:=\sup\{|\l|\mid \l\textrm{ is an eigenvalue of }A \}=\sup\{|<vA,v>| \mid |v|\leq 1\}$; recall also that $|A|\leq n|A|_*$).
For  $\nu=(\bar \a+2\e)^{-1}$, relation~\eqref{cil3} entails
\begin{equation}\label{T9}
-2(\bar\a+\e)\left(\begin{array}{cc} I&0\\0&I\end{array}\right)\le
\left(\begin{array}{cc} X&0\\0&-Y\end{array}\right)\le3\bar\a\left(\begin{array}{cc} I&-I\\-I&I\end{array}\right) + 2\e\left(\begin{array}{cc} I&0\\0&I\end{array}\right).
\end{equation}
%
%
From this inequality, one can deduce that, for every $(\th,\z)\in \Theta\times Z$, there holds
\begin{multline}\label{Ishii}
\tr\left(A_{i_0}^{\tz,1}(\t_0,x_0,p_{x_0})X\right)- \tr\left(A_{i_0}^{\tz, 2}(\t_0,y_0,p_{y_0})Y\right)\leq 3\bar\a \left|a_{i_0}^{\tz,1}(\t_0,x_0,p_{x_0})-a_{i_0}^{\tz,2}(\t_0,y_0,p_{y_0})\right|^2\\+
2\e (|a_{i_0}^{\tz,1}(\t_0,x_0,p_{x_0})|^2+|a_{i_0}^{\tz,2}(\t_0,y_0,p_{y_0})|^2)
\end{multline}
In order to prove this inequality, we shall use the arguments by Ishii~\cite{Is89}.
Multiplying the latter inequality in \eqref{T9} by the matrix
$$
\left(\begin{array}{lc}
a_{i_0}^{\tz,1}(\t_0,x_0,p_{x_0})a_{i_0}^{\tz,1}(\t_0,x_0,p_{x_0})^T&a_{i_0}^{\tz,2}(\t_0,y_0,p_{y_0})a_{i_0}^{\tz,1}(x_0,p_{x_0})^T\\[8pt]
a_{i_0}^{\tz,1}(\t_0,x_0,p_{x_0})a_{i_0}^{\tz,2}(\t_0,y_0,p_{y_0})^T&a_{i_0}^{\tz,2}(\t_0,y_0,p_{y_0})a_{i_0}^{\tz,2}(\t_0,y_0,p_{y_0})^T
\end{array}\right)
$$
(which is symmetric and nonnegative definite) and evaluating the trace, we obtain
\begin{align*}
&\tr\left(A_{i_0}^{\tz,1}(\t_0,x_0,p_{x_0}) X-A_{i_0}^{\tz,2}(\t_0,y_0,p_{y_0}) Y\right)\leq \\&\qquad
3\bar\a\tr\left[\left(a_{i_0}^{\tz,1}(\t_0,x_0,p_{x_0})-a_{i_0}^{\tz,2}(\t_0,y_0,p_{y_0})\right)\left(a_{i_0}^{\tz,1}(\t_0,x_0,p_{x_0})
-a_{i_0}^{\tz,2}(\t_0,y_0,p_{y_0})\right)^T\right]\\
&\qquad+2\e \tr\left(a_{i_0}^{\tz,1}(\t_0,x_0,p_{x_0})a_{i_0}^{\tz,1}(\t_0,x_0,p_{x_0})^T + a_{i_0}^{\tz,2}(\t_0,y_0,p_{y_0})a_{i_0}^{\tz,2}(\t_0,y_0,p_{y_0})^T\right)
\end{align*}
and therefore, by using our choice of $\e$ and of $\bar \a$, we get relation~\eqref{Ishii}.
Since $u^1$ is a subsolution to problem~\eqref{P}, the former relation in~\eqref{cil1}  and ~\eqref{cil2}
yield
\begin{equation}\label{T10}
\begin{split}
0&\geq a+\min\limits_{\z\in Z}\max\limits_{\th\in \Theta} \{
-\tr\left(A_{i_0}^{\tz,1}(\t_0,x_0,p_{x_0})X\right)+f_{i_0}^{\tz,1}(\t_0,x_0,u^1(\t_0,x_0),p_{x_0},X)\}\\
&\geq b +\frac{ \d \s}{t}+\min\limits_{\z\in Z}\max\limits_{\th\in \theta} \{
-\tr\left(A_{i_0}^{\tz,2}(\t_0,y_0,p_{y_0})Y\right)+f_{i_0}^{\tz,2}(\t_0,y_0,u^2(\t_0,y_0),p_{y_0},Y)\\
&\qquad+\tr\big(A_{i_0}^{\tz,2}(\t_0,y_0,p_{y_0})Y-A_{i_0}^{\tz,1}(\t_0,x_0,p_{x_0})X\big)+
f_{i_0}^{\tz,1}(\t_0,x_0,u^1(\t_0,x_0),p_{x_0},X)\\
 &\qquad-f_{i_0}^{\tz,2}(\t_0,y_0,u^2(\t_0,y_0),p_{y_0},Y) \}+\frac \a2 \bar \g e^{\bar \g\t_0}|x_0-y_0|^2.
\end{split}
\end{equation}
From \eqref{T9}, it follows that
\begin{equation}\label{cil4}
X\le Y+4\e I, \qquad |X|, |Y|\le n(3\bar \a+2\e);
\end{equation}
actually, in order to prove these estimates, it suffices to evaluate inequality \eqref{T9} on the vectors $(v,v)$, $(v,0)$ and $(0,v)$ respectively. Whence, assumption $(C2)$ ensures
\begin{equation}\label{cil4_1}
f_{i_0}^{\tz,2}(\t_0,y_0,u^2(\t_0,y_0),p_{y_0},Y)\le f_{i_0}^{\tz,2}(\t_0,y_0,u^2(\t_0,y_0),p_{y_0},X-4\e I).
\end{equation}
Moreover, by $\psi(\t_0,x_0,y_0,i_0)\ge \psi(\t_0,x_0,y_0,j)$ and by \eqref{T5},  we get respectively
\[u_{i_0}^1(\t_0,x_0 )-u_{i_0}^2(\t_0,y_0)\ge u_{j}^1(\t_0,x_0 )- u_{j}^2(\t_0,y_0)\qquad\text{ for any $j\in I$}
\]
and
\[u_{i_0}^1(\t_0,x_0 )-u_{i_0}^2(\t_0,y_0)\ge  0.\]
 Hence by ($C4$) and recalling that  $\g=\min(\g_R^1,\g_R^2)=0$, we get
\begin{equation}\label{T11}
 f_{i_0}^{\tz,2}(\t_0,y_0,u^2(\t_0,y_0),p_{y_0},X-4\e I)\le  f_{i_0}^{\tz,2}(\t_0,y_0,u^1(\t_0,x_0),p_{y_0},X-4\e I).
\end{equation}
By \eqref{T10}, \eqref{cil4_1} and \eqref{T11} we get
\begin{align*}
0&\geq b+\frac{ \d \s}{t}+\min\limits_{\z\in Z}\max\limits_{\th\in \Theta} \{
-\tr\left(A_{i_0}^{\tz,2}(\t_0,y_0,p_{y_0})Y\right)+f_{i_0}^{\tz,2}(\t_0,y_0,u^2(\t_0,y_0),p_{y_0},Y)\\
&\qquad+\tr\big(A_{i_0}^{\tz,2}(\t_0,y_0,p_{y_0})Y-A_{i_0}^{\tz,1}(\t_0,x_0,p_{x_0})X\big)+
f_{i_0}^{\tz,1}(\t_0,x_0,u^1(\t_0,x_0),p_{x_0},X)\\
&\qquad-f_{i_0}^{\tz,2}(\t_0,y_0,u^1(\t_0,x_0),p_{y_0},X-4\e I)\}+\frac \a2 \bar \g e^{\bar \g\t_0}|x_0-y_0|^2\\
&\ge b+\frac{ \d \s}{t}+\min\limits_{\z\in Z}\max\limits_{\th\in \Theta} \{
-\tr\left(A_{i_0}^{\tz,2}(\t_0,y_0,p_{y_0})Y\right)+f_{i_0}^{\tz,2}(\t_0,y_0,u^2(\t_0,y_0),p_{y_0},Y)\}
\\&\qquad+\min_{\th,\z}\{\tr\big(A_{i_0}^{\tz,2}(\t_0,y_0,p_{y_0})Y-A_{i_0}^{\tz,1}(\t_0,x_0,p_{x_0})X\big) +f_{i_0}^{\tz,1}(\t_0,x_0,u^1(\t_0,x_0),p_{x_0},X)\\
&\qquad-f_{i_0}^{\tz,2}(\t_0,y_0,u^1(\t_0,x_0),p_{y_0},X-4\e I)\}+\frac \a2 \bar \g e^{\bar \g\t_0}|x_0-y_0|^2.
\end{align*}
Hence, since $u^2$ is a supersolution, we get
\begin{equation}\label{T10bis}
\begin{split}
\frac{ \d \s}{t}&\le \max_{\th,\z}\{-\tr\big(A_{i_0}^{\tz,2}(\t_0,y_0,p_{y_0})Y-A_{i_0}^{\tz,1}(\t_0,x_0,p_{x_0})X\big)-
f_{i_0}^{\tz,1}(\t_0,x_0,u^1(\t_0,x_0),p_{x_0},X)\\
&\qquad+f_{i_0}^{\tz,2}(\t_0,y_0,u^1(\t_0,x_0),p_{y_0},X-4\e I)\}+\frac \a2 \bar \g e^{\bar \g\t_0}|x_0-y_0|^2\\
&\le \max\limits_{\th,\z}\Big\{3\bar \a\Big |a_{i_0}^{\tz,1}(\t_0,x_0,p_{x_0})-a_{i_0}^{\tz,2}(\t_0,y_0,p_{y_0})\Big|^2 \\
&\qquad+ f_{i_0}^{\tz,2}(\t_0,y_0,u^1(\t_0,x_0),p_{y_0},X-4\e I) -f_{i_0}^{\tz,1}(\t_0,x_0,u^1(\t_0,x_0),p_{x_0},X)\\
&\qquad+2\e (|a_{i_0}^{\tz,1}(\t_0,x_0,p_{x_0})|^2+ |a_{i_0}^{\tz,2}(\t_0,y_0,p_{y_0}) |^2) \Big\}+\frac \a2 \bar \g e^{\bar \g\t_0}|x_0-y_0|^2
\end{split}
\end{equation}
where the last inequality is due to  ~\eqref{Ishii}.\\
Set $p:=\a e^{\bar \g \t_0}(x_0-y_0)$, $p^x:=\e x_0$, $p^y:=\e y_0$ and observe that $p_{x_0}=p+p^x$, $p_{y_0}=p-p^y$.
We define
\begin{multline}\label{T12}
    F^{\a,\e}_t:=\{(\t,x,y,i,r,p,p^x,p^y,X,\th,\z):\, X=X_1+X_2,\,(\t,x,y,i,r,p,X_1,\th,\z)\in D^{\a}_{0t}\\
   |X_2|\leq 2n\e,\, \e^{1/2}|x|, \e^{1/2}|y|\le m(\e),|p^x|,|p^y|\le (2R\e)^{1/2}\}
\end{multline}
From \eqref{T10bis}, we have
\begin{align*}
\frac{ \d \s}{t} &\le \sup_{F^{\a,\e}_{t}}\Big\{ (3\bar \a+2\e)\big|a_{i}^{\tz,1}(\t,x,p+p^x)-a_{i}^{\tz,2}(\t,y,p-p^y)\big|^2\\
&\qquad+f_{i}^{\tz,2}(\t,y,r,p-p^y,X-4\e I) -  f_{i}^{\tz,1}(\t ,x,r ,p+p^x,X)  \\
&\qquad-\frac{\a}{2}\bar\g e^{\bar \g\t}|x-y|^2 +2\e (|a_{i}^{\tz,1}(\t,x,p+p^x)|^2+|a_{i}^{\tz,2}(\t ,y,p-p^y)|^2)\Big\}^+.
   \end{align*}
By definition of $F^{\a,\e}_t$ and $(C1)$ and $(C3)$, we get that there exists a modulus of continuity $\o$ such that
\begin{equation*}
\begin{split}
   \frac{ \d \s}{t}&\le \sup_{F^{\a,\e}_{t}}\Big\{ f_{i}^{\tz,2}(\t,y,r,p ,X) -  f_{i}^{\tz,1}(\t ,x ,r ,p,X)+ 3\bar \a\big|a_{i}^{\tz,1}(\t,x,p)-a_{i}^{\tz,2}(\t,y,p)\big|^2\\
&\qquad -\frac{\a}{2}\bar\g e^{\bar \g\t}|x-y|^2+\o(|p^x|+|p^y|+\e) +2\e(|a_{i}^{\tz,1}(\t,x,p+p^x)|^2
  +|a_{i}^{\tz,2}(\t ,y,p-p^y)|^2)\Big\}^+.
  \end{split}
   \end{equation*}
If $(\t,x,y,i)\in E^{\a}_t$, by definition of $\s$, we get
\[u^1_i(\t,x)-u^2_i(\t,y)-\frac{\a}{2} e^{\bar \g\t}|x-y|^2\le \s+\s_0+\e\left\{\frac 1 {t-\t} +\frac 1 2   (|x|^2+|y|^2 )  \right\}\]
By the last two inequalities we get
\begin{equation}\label{T13}
\begin{split}
& u^1_i(\t,x)-u^2_i(\t,y)-\frac{\a}{2} e^{\bar \g\t}|x-y|^2\\
&\quad\le \s_0+\frac t \d\sup_{F^{\a,\e}_{t}} \big\{f_{i}^{\tz,2}(\t,y,r,p ,X) -  f_{i}^{\tz,1}(\t ,x ,r ,p,X)
+  3\bar \a\big|a_{i}^{\tz,1}(\t,x,p)-a_{i}^{\tz,2}(\t,y,p)\big|^2\\
&\qquad -\frac{\a}{2}\bar\g e^{\bar \g\t}|x-y|^2+\o(|p^x|+|p^y|+\e)+2\e(|a_{i}^{\tz,1}(\t,x,p+p^x)|^2+|a_{i}^{\tz,2}(\t ,y,p-p^y)|^2)  \big\}^+\\
&\qquad +\e\left\{\frac 1 {t-\t} +\frac 1 2   (|x|^2+|y|^2 )  \right\}
  \end{split}
\end{equation}
Observe that, by $(C3)$ and the definition of $F^{\a,\e}_{t}$, we have
\[\e\left(|a_{i}^{\tz,1}(\t,x,p+p^x)|^2+|a_{i}^{\tz,2}(\t,x,p-p^y)|^2\right)\le  C\e(1+|x|^2)\le C\e(1+\frac{m(\e)^2}{\e}). \]
Then sending $\e\to 0$ (note that, since $\|u_1\|=\|u_1\| \wedge \|u_2\|$, by \eqref{T12}, as $\e \to 0$, $F^{\a,\e}_t$ converges to $D^\a_{0 t}$) and then $\d\to 1^-$ in \eqref{T13}, we get the estimate.\par
The general case $\g\neq 0$ can be proved following the argument of the corresponding result in \cite[Thm 3.1]{JK1}.
\end{Proof}

\begin{Remark}
This result can be generalized (using the same proof) to the case of a bounded above subsolution $u^1$ and a bounded below supersolution $u^2$. Actually, in order to treat this case, it suffices to set $R:=\max\{\|u^{1+}\|, \|u^{2-}\|\}$ and to replace ``$\min\{\|u^{1}\|, \|u^{2}\|\}$'' with ``$\min\{\|u^{1+}\|, \|u^{2-}\|\}$'' in the definition of $D^\a_{\g t}$.
\end{Remark}

\begin{Remark}
Theorem \ref{thm:dipcnt} can be improved when either $u^1$ or $u^2$ satisfies additional regularity properties.  For instance, when they are both continuous functions with modulus of continuity $\o_1$ and $\o_2$, respectively, the result of Theorem \ref{thm:dipcnt} holds with
$\Delta_\a$ in \eqref{Deltaalpha}  defined by
\[\Delta^\a:=\Big\{(x,y)\in\R^n\times \R^n:\,\a |x-y|^2-\o_1(|x-y|)-\o_2(|x-y|)\le 0 \Big\}.\]
When either $u^1$ or $u^2$ belongs to $C^1(Q_T)$, then $\Delta_\a$ can be defined as
\[\Delta^\a:=\Big\{(x,y)\in\R^n\times \R^n:\,|x-y|\le n\min([u^1]_1,[u^2]_1)\a^{-1} \Big\}.\]
With the previous definitions of $\D_\a$, the proof of Theorem \ref{thm:dipcnt} can be easily adapted by using Lemma \ref{Lemma2}.(ii) and (iii)
(see \cite{JK1} for more details).
\end{Remark}

\subsection{Systems arising in control theory}\label{oc_pb}
Weakly coupled systems are the dynamic programming equations of optimal control problem of Markov process with random switching (see \cite{fz98}) and arise  in many areas as in connection with the optimal control of hybrid systems (\cite{BMP,DHM, GAM}).
Consider  the control problem with dynamics
\begin{equation}\label{dyn}
\begin{split}
&dX (s)=  b^{\theta_s,\zeta_s}_{\nu_s}\Big(s, X(s)\Big)ds+a^{\theta_s,\zeta_s}_{\nu_s}\Big(s, X(s)\Big)dW_s,\qquad s\in [t,T]\\
&X (t)=x
\end{split}
\end{equation}
where  $W_t $ is a standard Brownian motion, $\th_t$, $\zeta_t$ are  the controls and
 $\nu_t$ is a continuous time random process with state space $\{1,\dots, m\}$ for which
\begin{equation}\label{nu}
\mathbb{P} \{\nu_{t+\Delta t}=j\,|\, \nu_t=i,\,X_t=x\}=c^{\th,\zeta}_{ij}(t,x)\Delta t+O(\Delta t)
\end{equation}
for $\Delta t\to 0$, $i,j\in I$, $i\neq j$.
Let  $v=(v_1,\dots,v_m)$ be the  value function defined by
\begin{equation}\label{Value}
    v_i(x,t)=\inf_{\th \in \mathcal{T}}\sup_{\zeta \in \mathcal{Z}} \E_{x,i}\left\{\int_t^T l_i^{\theta_s,\zeta_s}\Big(s, X (s)\Big)ds+ u_{0,i}(X (T))\right\}
    \qquad (t,x)\in [0,T)\times \R^n,\ \ i\in I.
\end{equation}
where $\mathcal{T}$ stands for the set of admissible strategies of the first player (namely, non-anticipating  maps $\th:\mathcal{Z}\to \mathcal{T}$, see     \cite{FSu}).
Then the function $u(x,T-t):=v(x,t)$ is \emph{formally} the solution of
\eqref{P} with initial datum $u_i(x,0)=u_{0,i}(x)$ where  the operators $H_i$ are defined by
\begin{equation}\label{hjbi}
H_i(t,x,r,p,X)=\min\limits_{\z\in Z}\max\limits_{\th\in \Th}
\left\{
-\tr\left(A_i^\tz(t,x)X\right)+b^\tz_i(t,x)\cdot p+l^\tz_i(t,x)+\sum_{j\in I}d_{ij}^\tz(t,x)r_j\right\}
\end{equation}
and $A_i^\tz(t,x)=a^{\tz}(t,x)a^{\tz}(t,x)^T$, $d^\tz_{ij}=-c^\tz_{ij}$ for $j\neq i$, $d^\tz_{ii}=-\sum_j c^\tz_{ij}$.
Besides assumptions $(C0)$-$(C6)$, we require the following assumptions
\begin{itemize}
\item the coefficients $a^\tz_i$ and  $b^\tz_i$ are Lipschitz continuous in $x$ uniformly in $t$ and in $(\th,\z)$; namely, for $h=a^\tz_i,b^\tz_i$, there holds
\begin{equation*}
|h(t,x)-h(t,y)|\leq L_h|x-y| \qquad \forall x,y,t,\th,\z,i;
\end{equation*}
\item the coefficients $l^\tz_i$ and $d^\tz_{ij}$ are H\"older continuous in $x$ uniformly in $t$ and in $(\th,\z)$; namely, for $ h=l^\tz_i,d^\tz_{ij}$, there holds
\begin{equation*}
|h^\tz_i(t,x)-h^\tz_i(t,y)|\leq L_{h }|x-y|^\mu \qquad \forall x,y,t,\th,\z,i.
\end{equation*}
\end{itemize}

\begin{Theorem}\label{cde_hjbi}
Let $H^1$ and $H^2$ be two operators of the form \eqref{hjbi} which fulfill the above assumptions.
Let $u^1$ and $u^2$ be respectively a bounded subsolution  to problem~\eqref{pbC} with $H=H^1$ and $u_0=u^1_0\in C^\mu(\R^n)$ and a bounded supersolution to problem~\eqref{pbC} with $H=H^2$ and $u_0=u^2_0\in C^\mu(\R^n)$. Set $R:=\max (\|u^1\|,\|u^2\|)$ and $\g:=\min\{\g^1_R,\g^2_R\}$.
Then there exists a constant $K>0$ (depending only on $T$, $R$ and on the constants entering in our assumptions) such that, for every $0\le t\le T$, there holds
\begin{multline*}
e^{\g t}\|u^1(t,\cdot)-u^2(t,\cdot)\|
\leq \|u^1_0-u^2_0\| +K t \sup_{i\in I,(\t,x)\in Q_t,\th,\z}\left[|\ell^{\tz,1}_{i}-\ell^{\tz,2}_{i}|+|d^\tz_1-d^\tz_2|\right]\\
+K t^{\mu/2}\sup_{i\in I,(\t,x)\in Q_t,\th,\z}\left[|b^{\tz,1}_{i}-b^{\tz,2}_{i}|^\mu +|a^{\tz,1}_{i}-a^{\tz,2}_{i}|^\mu\right]
\end{multline*}
where $d^{\tz}_k$ is the matrix $(d^{\tz,k}_{ij})_{i,j\in I}$ for $k=1,2$.
\end{Theorem}

\begin{Remark}
If $u^1$ and $u^2$ are both solutions and $\sup_{\zeta,\theta,t,x,i}|l^{\tz,k}_{i}(t,x)|<+\infty$ ($k=1,2$), then Theorem~\ref{bound} below guarantees that $u^1$ and $u^2$ are bounded and it also provides an estimate of $R$.
\end{Remark}

\begin{Proof}
By the arguments of~\cite[Thm 3.2 and 4.1]{JK1}, this result is a consequence of Theorem~\ref{thm:dipcnt} and of the regularity of the coefficients and of the initial data.

For $C_\mu:=\min\{[u^1_0]_\mu,[u^2_0]_\mu\}$, we have
\begin{equation*}
e^{\g t}\|u^1(t,\cdot)-u^2(t,\cdot)\|\leq
\sup_{E_t^\a}\Big(e^{\g \t}\big(u^1_i(\t,x)-u^2_i(\t,y)\big) -\frac \a 2e^{\bar  \g \t}|x-y|^2\big)\Big)
\end{equation*}
\begin{equation*}
\sup_{E_0^\a}  \big(u^1_{0,i}(x)-u^2_{0,i}(y) -\frac \a 2|x-y|^2\big)^+ \leq
\|u^1_0-u^2_0\|+ \a^{-\mu/(2-\mu)}C_\mu^{2/(2-\mu)}.
\end{equation*}
By these inequalities, for
\begin{equation*}
f^{\tz,k}_i=b^{\tz,k}_i(t,x)\cdot p+l^{\tz,k}_i(t,x)+\sum_{j\in I}d^{\tz,k}_{ij}(t,x)r_j, \qquad a^{\tz,k}=a^{\tz,k}(t,x)\qquad (k=1,2),
\end{equation*}
Theorem~\ref{thm:dipcnt} yields
\begin{align*}
&e^{\g t}\|u^1(t,\cdot)-u^2(t,\cdot)\| \\
&\quad\leq \|u^1_0-u^2_0\|+ \a^{-\mu/(2-\mu)}C_\mu^{2/(2-\mu)}
+ t\sup_{D^\a_{\g t}} \Big\{e^{\g \t}\Big[|b^{\tz,1}_{i}(\t,x)-b^{\tz,2}_{i}(\t,x)||p|
\\&\qquad+|l^{\tz,1}_{i}(\t,x)-l^{\tz,2}_{i}(\t,x)|
+R\sum_{j\in I}|d^{\tz,1}_{ij}(\t,x)-d^{\tz,2}_{ij}(\t,x)|
+L_b|x-y||p| 
\\&\qquad+(L_l+mL_d R)|x-y|^\mu\Big]
+ 3\a e^{\bar \g\t}\left(|a^{\tz,1}_{i}(\t,x)-a^{\tz,2}_{i}(\t,x)|+L_a|x-y|\right)^2
-\frac \a 2\bar\g e^{\bar \g \t}|x-y|^2\Big\}
\\&
\quad\leq \|u^1_0-u^2_0\|+ \a^{-\mu/(2-\mu)}C_\mu^{2/(2-\mu)}
+ t\sup_{D^\a_{\g t}} \Big\{e^{\g \t}\Big[
|l^{\tz,1}_{i}(\t,x)-l^{\tz,2}_{i}(\t,x)|
\\&\qquad +R\sum_{j\in I}|d^{\tz,1}_{ij}(\t,x)-d^{\tz,2}_{ij}(\t,x)|\Big]
+\a e^{\bar \g \t}|b^{\tz,1}_{i}(\t,x)-b^{\tz,2}_{i}(\t,x)|^2 
\\&\qquad
+e^{\g\t}(L_l+mL_d R)|x-y|^\mu
+6\a e^{\bar \g\t}|a^{\tz,1}_{i}(\t,x)-a^{\tz,2}_{i}(\t,x)|^2
\\&\qquad
+\a e^{\bar \g \t}|x-y|^2\left(1+L_b+6L_a^2-\frac{\bar \g}2\right)\Big\}
\end{align*}
where the last inequality is due to the Young one and to the choice $p=\a (x-y)e^{(\bar \g -\g)\t}$.

We choose $\bar \g$ sufficiently large such that
\begin{equation*}
1+L_b+ 6L_a^2-\bar\g/2=-1.
\end{equation*}
Furthermore, by standard calculus, we get
\begin{equation*}
(L_l+mL_d R)e^{\g \t}|x-y|^\mu-\a e^{\bar \g \t}|x-y|^2\leq
K_1\a^{-\mu/(2-\mu)}
\end{equation*}
where $K_1$ is a constant depending only on $L_l$, $L_d$,$R$, $\g$, $\bar \g$ and $T$.
Taking into account the last three inequalities, for $K_2:=C_\mu^{2/(2-\mu)/2}+tK_1$ we obtain
\begin{align*}
&e^{\g t}\|u^1(t,\cdot)-u^2(t,\cdot)\|\\
&\quad
\leq \|u^1_0-u^2_0\|+ t\sup_{D^\a_{\g t}} \Big\{e^{\g \t}\Big[|l^{\tz,1}_{i}(\t,x)-l^{\tz,2}_{i}(\t,x)|+R\sum_{j\in I}|d^{\tz,1}_{ij}(\t,x)-d^{\tz,2}_{ij}(\t,x)|
\Big]\Big\}\\
&\qquad+K_2 \a^{-\mu/(2-\mu)}+\a t \sup_{D^\a_{\g t}} \Big\{e^{\bar \g \t}\left[|b^{\tz,1}_{i}(\t,x)-b^{\tz,2}_{i}(\t,x)|^2 +6|a^{\tz,1}_{i}(\t,x)-a^{\tz,2}_{i}(\t,x)|^2\right]\Big\}
\end{align*}
Minimizing the right-hand side by an adequate choice of $\a$, we have
\begin{align*}
&e^{\g t}\|u^1(t,\cdot)-u^2(t,\cdot)\|\\
&\quad \leq \|u^1_0-u^2_0\|+ t\sup_{D^\a_{\g t}} \Big\{e^{\g \t}\Big[|l^{\tz,1}_{i}(\t,x)-l^{\tz,2}_{i}(\t,x)|+R\sum_{j\in I}|d^{\tz,1}_{ij}(\t,x)-d^{\tz,2}_{ij}(\t,x)|
\Big]\Big\}\\
&\qquad+K_3 t^{\mu/2}\left( \sup_{D^\a_{\g t}} \Big\{e^{\bar \g\t}\left[|b^{\tz,1}_{i}(\t,x)-b^{\tz,2}_{i}(\t,x)|^2 +|a^{\tz,1}_{i}(\t,x)-a^{\tz,2}_{i}(\t,x)|^2\right]\Big\}\right)^{\mu/2}.
\end{align*}
where $K_3$ is a constant depending only on $T$, and on the constants entering in our assumptions.
Finally, by the Young inequality, one can easily accomplish the proof.
\end{Proof}
\begin{Remark}
  By the proof above, this result is still true when, for each $i\in I$, either $u^1_{0,i}$ or $u^2_{0,i}$ belongs to $C^\mu(\R^n)$ and the assumptions on the regularity are fulfilled by either $a^{\tz,1}_i$ or $a^{\tz,2}_i$, either $b^{\tz,1}_i$ or $b^{\tz,2}_i$, either $l^{\tz,1}_i$ or $l^{\tz,2}_i$.
\end{Remark}

%
%

\section{Regularity estimates and vanishing viscosity}\label{reg_vv}
In this Section we collect some applications of Theorem~\ref{thm:dipcnt}:
the first part is devoted to establish a regularity estimate for the solution to system~\eqref{P} provided that the initial condition and the coefficients are H\"older continuous. In the second part we prove an estimate of the vanishing viscosity approximation.
\subsection{Regularity estimates}
In this section, we address the Cauchy problem
\begin{equation}\label{pbC}
\left\{\begin{aligned}
&\partial_t u_i+H_i(t,x,u,Du_i,D^2u_i)=0 &&\qquad \textrm{in }Q_T\\
&u_i(0,x)=u_{0i}(x) &&\qquad\textrm{on } \R^n, \, i\in I
\end{aligned}\right.
\end{equation}
with $H_i$ of the form \eqref{minmax} and we establish two results for the solution $u$: an $L^\infty$-estimate and the H\"older continuity.
\begin{Theorem}\label{bound}
Assume conditions~($C0$)-($C4$) and ($C7$).
For $u_{0i}$ continuous and bounded ($i\in I$), let $u$ be a bounded solution to problem~\eqref{pbC}. Then, for $\g:=\g_{\|u\|}$ (the constant~$\g_R$ is introduced in (C4)) there holds:
\begin{equation*}
\|u(t,\cdot)\| \leq e^{-\g t}\|u_0\|+ t e^{\g t} C^f.
\end{equation*}
\end{Theorem}
\begin{Proof}
Assume $\g_{\|u\|}=0$ in ($C4$). We shall proceed following the same arguments as those of Theorem~\ref{thm:dipcnt} with $u_1\equiv 0$ and $u_2=u$ (clearly, $u_1$ is the solution to \eqref{P} with zero coefficients). Relations \eqref{cil1} and \eqref{cil2} guarantee
\begin{eqnarray*}
0&\leq& b+\min\limits_{\z\in Z}\max\limits_{\th\in \Theta} \{
-\tr\left(A_{i_0}^{\tz}(\t_0,y_0,p_{y_0})Y\right)+f_{i_0}^{\tz}(\t_0,y_0,u(\t_0,y_0),p_{y_0},Y) \}\\
&\leq&
 a-\frac{\d \s}t +\min\limits_{\z\in Z}\max\limits_{\th\in\Theta} \{
-\tr\left(A_{i_0}^{\tz}(\t_0,y_0,p_{y_0})Y\right)+f_{i_0}^{\tz}(\t_0,y_0,u(\t_0,y_0),p_{y_0},Y) \}
\end{eqnarray*}
where $p_{y_0}=\a e^{\bar\g \t_0}(x_0-y_0)-\e y_0$. We observe that $(a,p_{x_0},X)\in \bar{\cal P}^{2,+}0$ iff $a=0$, $p_{x_0}=0$ and $X\geq 0$; hence, by~\eqref{cil4}, we get $Y\geq X-4\e I\geq -4\e I$. Therefore, the above estimate entails
\begin{eqnarray*}
\s &\leq& t \d^{-1}\min\limits_{\z\in Z}\max\limits_{\th\in\Theta} \{
4\e\tr\left(A_{i_0}^{\tz}(\t_0,y_0,p_{y_0})\right)+f_{i_0}^{\tz}(\t_0,y_0,u(\t_0,y_0),p_{y_0},-4\e I) \}\\
&\leq & t \d^{-1}\min\limits_{\z\in Z}\max\limits_{\th\in\Theta} \{
4\e\tr\left(A_{i_0}^{\tz}(\t_0,y_0,\textrm{const.}\e^{1/2})\right)+f_{i_0}^{\tz}(\t_0,y_0,0,\textrm{const.}\e^{1/2},-4\e I) \}
\end{eqnarray*}
where the last inequality is due to the same arguments as in \eqref{T11} and to Lemma \ref{Lemma2}-(i) and -(iii). Observe that assumption ($C3$) and estimate \eqref{T8} ensure
\begin{equation*}
\e \tr\left(A_{i_0}^{\tz}(\t_0,y_0,\textrm{const.} \e^{1/2})\right)\leq \textrm{const.} m(\e).
\end{equation*}
 Letting $\e\to 0$ and $\d\to 1$, we obtain
\begin{equation*}
-\|u_0\|-\inf\limits_{x\in\Rn,i\in I}\{u_i(t,\cdot)\}\leq \s\leq t C^f;
\end{equation*}
namely, one side of the statement is established. Reversing the role of $u$ and 0, one can easily deduce the other inequality of the statement. The case $\g_{\|u\|}\neq 0$ will follow as for Theorem~\ref{thm:dipcnt}.
\end{Proof}

\begin{Theorem}\label{regularity}
Assume  $(C0)$-$(C6)$ and $u_{0}\in C^\mu(\R^n)$ for some $\mu\in(0,1]$. Then any bounded solution~$u$ to problem~\eqref{pbC} is H\"older continuous in~$x$ and, for some positive constant~$K$, it fulfills
$$
[u(t,\cdot)]_\mu\leq K e^{\bar \g t} \left([u_0]_\mu+t^{1-\mu/2}e^{\g^+t}C_{f,\|u\|}\right)
$$
where $\g^+:=\max\{0,\g_{\|u\|}\}$ and $\bar \g:=2(C_{f,\|u\|}+ 3C_a^2+1)+\g^+$ (the constants $\g_R$, $C_{f,R}$ and $C_a$ are those introduced respectively in $(C4)$, $(C5)$ and $(C6)$).
\end{Theorem}

\begin{Proof}
This proof relies on the arguments of \cite[Thm 3.3-(b)]{JK1}: the application of Theorem~\ref{thm:dipcnt} with $f^{\tz,1}_{i}=f^{\tz,2}_{i}$, $a^{\tz,1}_{i}=a^{\tz,2}_{i}$ and $u^1=u^2$ with a careful estimates of the two sides. For the sake of completeness, let us sketch them.
We observe that
\begin{align*}
    \sup_{i\in I,\, x\in\R^n}\Big(e^{\gamma t}(u^1_i(t,x)-u^2_i(t,x))\Big)\le  \sup_{E_t^\a}\Big(e^{\g \t}\big(u^1_i(\t,x)-u^2_i(\t,y)\big) -\frac \a 2e^{\bar  \g \t}|x-y|^2\big)\Big)
\end{align*}
and
\begin{equation*}
\sup_{E_0^\a}  \big(u^1_i(0,x)-u^2_i(0,y) -\frac \a 2|x-y|^2\big)^+ \leq [u_0]_\mu |x-y|^\mu-\frac \a 2|x-y|^2\leq
2[u_0]_\mu^{\frac {2}{2-\mu}}\a^{\frac {-\mu}{2-\mu}}
\end{equation*}
where the last inequality is due to the Young inequality with exponents $2/\mu$ and $2/(2-\mu)$.
Moreover, by conditions $(C5)$ and $(C6$) (recall $p=\a(x-y)e^{(\bar\g-\g)\t}$) and by our choice of $\bar \g$, we have
\begin{align*}
&e^{\g \t}[f^{\tz,2}_{i}(\t,y,r,p,X)-f^{\tz,1}_{i}(\t,x,r,p,X)]+ 3\a e^{\bar \g\t}|a^\tz_{1,i}(\t,x,p)-a^\tz_{2,i}(\t,y,p)|^2-\frac \a 2\bar\g e^{\bar \g \t}|x-y|^2\\
&\quad\leq  e^{\g \t}C_{f,\|u\|}|x-y|^\mu +\a e^{\bar \g \t}|x-y|^2\left(C_{f,\|u\|}+ 3C_a^2-\frac{\bar \g}2 \right)\\
&\quad \leq  e^{\g \t}\left[C_{f,\|u\|}|x-y|^\mu -\a |x-y|^2\right]\\
&\quad\leq
K_1 e^{\g t}(C_{f,\|u\|})^{2/(2-\mu)}\a^{-\mu/(2-\mu)}
\end{align*}
where last inequality is due to standard calculus and $K_1$ is a constant depending only on $\mu$. Therefore, taking into account the last two inequalities, Theorem~\ref{thm:dipcnt} entails
$$
e^{\g \t}\big(u^1_i(\t,x)-u^2_i(\t,y)\big) \leq \left[2[u_0]_\mu^{\frac {2}{2-\mu}}+K_1 e^{\g t}(C_{f,\|u\|})^{2/(2-\mu)}\right]\a^{-\mu/(2-\mu)}+\frac \a 2e^{\bar  \g \t}|x-y|^2
$$
and the statement follows by a suitable choice of $\a$ (see \cite[Thm 3.3-(b)]{JK1} for detailed calculations).
\end{Proof}
\subsection{Vanishing viscosity}
We consider the viscous approximation to \eqref{P}
\begin{equation}\label{VV}
\partial_t u_i^\e+H_i(t,x,u,Du^\e_i,D^2u^\e_i)=\e \D u^\e_i \qquad\textrm{in }Q_T,\, i\in I
\end{equation}
where $H_i$ is as in \eqref{minmax}.  In the next proposition
we establish an estimate on the  rate of convergence of $u^\e$ to $u$.
\begin{Proposition}\label{vv}
Assume $(C0)$-$(C7)$ and that, for any $\e>0$, there exists a bounded solution $u^\e$ to \eqref{VV}. Then there exists a solution $u\in C^\mu(Q_T)$ to  \eqref{P}-\eqref{minmax} and
\[
    \|u( t,\cdot)-u^\e(t,\cdot)\|\le C\big(\|u(0,\cdot)-u^\e(0,\cdot)\|+ \e^{\mu/2}\big) \qquad t\in [0,T]
\]
where $C$ is independent of $\e$.
\end{Proposition}
\begin{Proof}
The existence of the solution $u$ to \eqref{P} and the local uniform convergence of the sequence $u^\e$ to $u$ can be obtained by employing the classical weak limit method introduced by Barles-Perthame, which can be easily adapted to systems.
Moreover by Theorem \ref{regularity}, the functions $u^\e$ and $u$ belong to $C^\mu(Q_T)$ for any $\e$.
The proof of the rate of convergence is based on the estimate in Theorem \ref{thm:dipcnt} applied to problem \eqref{P} and \eqref{VV} with
\begin{align*}
    &f^{\tz,1}_i(t,x,r,p,X)=f^\tz_i(t,x,r,p,X)\\
    &f^{\tz,2}_i(t,x,r,p,X)=f^{\tz}_i(t,x,r,p,X)-\e \tr(X)\\
    &A^{\tz,1}_i(t,x,p) = A^{\tz,2}_i(t,x,p)= A_i^\tz(t,x,p)
\end{align*}
Since it is very similar to the proof of the corresponding result in \cite{JK1},   we omit it.
\end{Proof}
\begin{Remark}
A similar estimate for  the vanishing viscosity approximation of weakly coupled systems has been recently proved in \cite{CGT} using
different techniques and stronger assumptions.
\end{Remark}

\subsubsection{Vanishing viscosity for a first order problem}
Let us establish a rate of convergence for the vanishing viscosity approximation of a first order system arising in optimal control problem. Being a straightforward application of Proposition \ref{vv}, the proof is omitted.
\begin{Proposition}
Assume the hypotheses of Section~\ref{oc_pb}. Let $u_\e$ and $u$ be the solution of
\[
\partial_t u^\e_i+\min\limits_{\z\in Z}\max\limits_{\th\in \Th}
\left\{
-\e \tr\left(A_i^\tz(t,x)D^2u_i^\e\right)+b^\tz_i(t,x)\cdot Du_i^\e+l^\tz_i(t,x)+\sum_{j\in I}d_{ij}^\tz(t,x)u^\e_j\right\}=0
\]
and respectively of
\[
\partial_t u_i +\min\limits_{\z\in Z}\max\limits_{\th\in \Th}
\left\{b^\tz_i(t,x)\cdot Du_i+l^\tz_i(t,x)+\sum_{j\in I}d_{ij}^\tz(t,x)u_j\right\}=0.
\]
Then
\[
    \|u( t,\cdot)-u^\e(t,\cdot)\|\le C\big(\|u(0,\cdot)-u^\e(0,\cdot)\|+ \e^{\mu/2}\big) \qquad t\in [0,T].
\]
\end{Proposition}

%
%
\section{Periodic  Homogenization of quasi-monotone systems}\label{homogenization}
In this section we study the periodic homogenization of the fully nonlinear systems
\begin{equation}\label{HJe}
\left\{\begin{aligned}
&\partial_t u_i^\e +H_i\left(x,\frac x\e, u^\e,Du_i^\e,D^2u_i^\e\right)=0 &&\qquad \textrm{in }Q_T\\
&u_i^\e(0,x)=u_{0i}(x) &&\qquad\textrm{on } \R^n, \, i\in I
\end{aligned}\right.
\end{equation}
where
\[H_i(x,y, r, p, X)= \min\limits_{\z\in Z}\max\limits_{\th\in \Th}
\left\{
- \tr\left(A_i^\tz(x,y)X\right)+f^\tz_i(x,y,r,p)\right\}.\]
For the sake of clarity, let us list the assumptions that will hold throughout this section.
\begin{itemize}
\item[($H0$)] The sets $\Th$ and $Z$ are two compact metric spaces.
\item[($H1$)] The functions $f^{\tz}_i$ are continuous and, for some constant $L_f$ and a modulus of continuity $\o$, they satisfy
 \begin{multline*}
|f_i^\tz(x_1,y_1,r_1,p_1)-f_i^\tz(x_1,y_1,r_1,p_1)|\le  L_f|(x_1,y_1)-(x_2,y_2)|(|p_1|\vee|p_2|+|r_1|\vee|r_2|)\\+\o(|(x_1,y_1)-(x_2,y_2)|)
    +L_f(|r_1-r_2|+|p_1-p_2|)
 \end{multline*}
for every $x_k,y_k,r_k,\th,\zeta, i$ ($k=1,2$).
Moreover, there exists a constant $C$ such that
 \begin{equation}\label{bd}
 |f_i^\tz(x,y,0,0)|\le C \qquad \forall x,y,r,\th,\zeta.
 \end{equation}
\item[($H2$)]  $A^\tz_i(x,y)=a^\tz_i(x,y) a^\tz_i(x,y)^T$ for some bounded, continuous matrix $a_i^\tz$ satisfying
$$
\left| a^\tz_i(x_1,y_1)-a^\tz_i(x_1,y_2)\right|\leq L_a|(x_1,y_1)-(x_2,y_2)|
\qquad\forall x_k,y_k,\th,\zeta, i\in I \quad (k=1,2).
$$
\item[($H3$)] $f^{\tz}_i(x,\cdot,r,p)$ and $a^\tz_i(x,\cdot)$ are $\mathbb Z^n$-periodic in $y$ for any $x,r,p,\th,\zeta,i$.
\item[($H4$)] The matrix $A^\tz_i$ is uniformly elliptic, namely, for some positive constant $\nu$ there holds
\begin{equation*}
a^\tz_i(x,y)\geq \nu I, \qquad \forall x,y,\th,\z,i.
\end{equation*}
\item[($H5$)]
There exists $\g\in\R $ such that if  $r$, $s\in\R^m$
and $\displaystyle r_j-s_j=\max_{k\in I}\{r_k-s_k\}\ge 0$, then
\[f^\tz_j(x,y,r,p)-f^\tz_j(x,y,s,p)\ge\gamma(r_j-s_j) \qquad \forall x,y,p, \th,\z.\]

\end{itemize}
We consider the {\it cell problem}:\\
For any fixed $i\in I$ and $(x,r,p,X)\in\R^n\times\R^m\times\R^n\times S^n$, find a constant $\Ho_i=\Ho_i(x,r,p,X)$ such that the equation
\begin{equation}\label{CP}
    H_i(x,y,r,p, X+D^2_yv(y))=\Ho_i,\qquad y\in\R^n
\end{equation}
admits a periodic solution $v_i= v_i(\cdot;x,r,p,X)$.

It is well known (see: \cite{Ev2, AL, AB8,m}) that there exists exactly one value $\Ho_i$ such that \eqref{CP} has a solution; moreover, $\Ho_i$ can be obtained as the (uniform) limit of $-\l v_{\l,i}$ as $\l\to 0$, where the {\it approximated corrector} $v_{\l,i}:=v_{\l,i}(y; x,r,p,X)$ is the solution to
\begin{equation}\label{eid}
   \l v_{\l,i}+H_i(x,y,r,p,X+D^2_y v_{\l,i})=0,\qquad y\in\R^n.
\end{equation}
We associate to each Hamiltonian $H_i$ the corresponding {\it effective Hamiltonian} $\Ho_i$.
Note that  at this level the index $i$ is fixed, hence the definition of the effective Hamiltonians does not involve any coupling among the equations. Nevertheless, in view of existence and uniqueness
results for the homogenized problem,  we  need to study
the regularity of the effective Hamiltonians in particular with respect to the variable $r\in \R^m$.

In the next proposition we collect some useful properties of the approximated correctors $v_{\l,i}$ and of the effective operators $H_i$.
\begin{Proposition}\label{Prp:erg}
The following properties hold:
\begin{itemize}
\item[$i$)] For any $i,x,r,p,X$, the  approximated equation~\eqref{eid} admits exactly one periodic continuous solution $v_{\l,i}$. Moreover, as $\l\to 0^+$, $\l v_{\l,i}$ and $(v_{\l,i}-v_{\l,i}(0))$ converge respectively to the ergodic constant~$-\Ho_i$ and to a solution~$v_i$ of~\eqref{CP} with $v_i(0)=0$.
\item[$ii$)]
For any $i\in I$, the effective Hamiltonian $\Ho_i$ is  continuous  in $(x,r,p,X)$ and
\begin{itemize}
\item[a)] For some constant $C_1>0$ and a modulus of continuity $\o_1$, there holds 
\begin{align*}
& |\Ho_i(x,r_1,p_1,X_1)- \Ho_i(x,r_2,p_2,X_2)|\leq C_1\left(|r_1-r_2|+|p_1-p_2|+|X_1-X_2|\right) ;\\
&   |\Ho_i(x_1,r,p,X)- \Ho_i(x_2,r,p,X)|\leq C_1(1+|p|+|r|+|X|)|x_1-x_2| +\o_1(|x_1-x_2|) ;\\
&  |\Ho_i(x,r,p,X)|\leq \max_{y,\th,\z}\left|-\tr(A^\tz_i(x,y)X)+f^\tz_i(x,y,r,p)\right|
\end{align*}
for every $x_k,p_k, r_k, X_k$ ($k=1,2$).
\item[b)] $\Ho_i$ is uniformly elliptic. Moreover, if $H_i$ is convex, then $\Ho_i$ is also convex.
\item[c)] $\{\Ho_i\}_{i\in I}$ is quasi-monotone, namely, it satisfies \eqref{QM}.
\end{itemize}
\end{itemize}
\end{Proposition}
\begin{Proof}
For statement ($i$), we refer to \cite{Ev2} (see also \cite{AB8} and \cite{AL}). The estimates in $(ii).a$     follow  by  the continuous dependance estimates in \cite[Thm 3.1]{m} (note that in the cell problem both $r$ and $p$ are fixed), while property $(ii).b$ is proved for example in \cite{AB8} and in \cite{Ev2}.
We finally prove that $\Ho_i$, $i\in I$, satisfy the quasi-monotonicity condition \eqref{QM}. Assume by contradiction that there exist  $r$, $s\in\R^m$ such that
$\displaystyle r_j-s_j=\mathop{\max}_{k\in I}\{r_k-s_k\}\ge 0$ and
\begin{eqnarray*}
\Ho_j(x, r,p, X)<\Ho_j(x ,s,p,X)
\end{eqnarray*}
for some $x\in \R^n$, $p\in\R^n$, $X\in S^n$. Let $u_r$ and $u_s$  be two periodic  solutions respectively of
\begin{align*}
    &H_j (x,y,r , p,X+ D^2 u_r)=\Ho_j(x,r,p,X)\qquad y\in \R^n,\\
     &H_j(x,y,s, p, X+D^2u_s)=\Ho_j(x,s,p,X)\qquad y\in \R^n.
\end{align*}
(these functions exist by point $(i)$).
Since $u_r$, $u_s$ are bounded, by adding a constant we can assume  w.l.o.g.  $u_r>u_s$ in $\R^n$.
Since
\begin{eqnarray*}
H_j (x,y,r ,p,X+ D^2 u_r)&=&\Ho_j(x,r,p,X)<\Ho_j(x,s,p,X)
=  H_j(x,y,s,p, X+D^2u_s)\\&\le& H_j(x,y,r,p, X+D^2u_s)
\end{eqnarray*}
(where the last inequality follows by $(H5)$), then    for $\l$ sufficiently small
\[
\l u_r+H_j (x,y,r ,Du_r,X+ D^2 u_r)\leq\l u_s+H_j(x,y,r,Du_s, X+D^2u_s) \qquad y\in\R^n.
\]
By the comparison principle for problem \eqref{eid}, we deduce $\l u_r\leq \l u_s$; as $\l\to 0^+$, we infer $\Ho_i(x,r,p,X)\geq \Ho_i(x,s,p,X)$ which gives the desired contradiction.
\end{Proof}
\begin{Proposition}\label{homo}
Let  $u_0\in BUC(\R^n)$. Then
\begin{itemize}
  \item For any $\e>0$ there exists a unique solution $u_\e\in BUC(Q_T)$
to \eqref{HJe}. Moreover $u_\e$ is bounded uniformly in $\e$.
  \item There exists a unique solution  $u\in BUC(Q_T)$ to the effective problem
\begin{equation}\label{hom}
\left\{\begin{aligned}
&\partial_t u_i +\Ho_i(x,u,Du_i, D^2u_i)=0 &&\qquad \textrm{in }Q_T\\
&u_i(0,x)=u_{0i}(x) &&\qquad\textrm{on } \R^n, \, i\in I
\end{aligned}\right.
\end{equation}
where the operators $\Ho_i$ are defined by the cell problem \eqref{CP}.
\end{itemize}
\end{Proposition}
\begin{Proof}
By routine adaptation of the arguments in \cite{IK91}, \eqref{HJe} and \eqref{hom} satisfy a comparison principle for sub and supersolution.\par

In order to prove the existence of the solution, we note that assumption $(H1)$ ensures $|f^\tz_i(x,y,r,0)|\leq C+L|r|$. We deduce that, for a constant $\tilde C$ sufficiently large, the functions $u^{\pm}(x,t)=\pm(\|u_0\|+e^{\tilde Ct},\dots,\|u_0\|+e^{\tilde Ct})$ are respectively a super- and a subsolution of \eqref{HJe}. Actually, by this inequality, we have
\begin{equation*}
\partial_t u^+_i +H_i\left(x,\frac x\e, u^+,Du_i^+,D^2u_i^+\right)=
\tilde C e^{\tilde Ct}+H_i\left(x,\frac x\e, u^+,0,0\right)\geq
(\tilde C-L) e^{\tilde Ct}-C-L\|u_0\|\geq0
\end{equation*}
provided that $\tilde C=L+1+C+L\|u_0\|$; hence $u^+$ is a supersolution. Being similar, the proof for $u^-$ is omitted.
By the Perron's method for system, see \cite{IK91}, it follows the existence
of a solution $u_\e\in BUC(Q_T)$ to \eqref{HJe} such that
\[-\|u_0\|-e^{\tilde CT}\le u^\e_i(t,x)\le\|u_0\|+e^{\tilde CT},\qquad (t,x)\in Q_T,\,i\in I.\]
The existence of a bounded solution to \eqref{hom} is proved in the same way.
\end{Proof}
\begin{Theorem}\label{thm:hom}
The solution $u^\e$ of \eqref{HJe} converges locally uniformly on $[0,T]\times \R^n$ to the solution $u\in BUC(Q_T)$  of \eqref{hom}.
\end{Theorem}
\begin{Proof}
By Proposition \ref{homo}  there exists a
continuous solution $u^\e $ of \eqref{HJe}
which is bounded independently of $\e.$ 
We follow the argument in \cite[Thm 3.5]{hi}. We introduce the half-relaxed
limits
\begin{eqnarray*}
\overline u (t,x) = \mathop{\rm lim\,sup}_{\e \to 0,  (t_\e,x_\e)   \to (t,x)} u^\e (t_\e,x_\e)
\ \ \ {\rm and} \ \ \
\underline u (t,x) = \mathop{\rm lim\,inf}_{\e \to 0,  (t_\e,x_\e)   \to (t,x) } u^\e (t_\e,x_\e).
\end{eqnarray*}
We first show that $\overline  u$ is a subsolution of the system \eqref{hom}.
We assume   there exists  $i\in I$ and $\phi\in C^2$ such that
$\overline{u}_i-\phi$ has a strict maximum point at some $(\tov,\xo)\in (0,T)\times \R^n$ with
$\overline u_i(\tov,\xo)=\phi(\tov,\xo).$ We assume wlog $i=1$ and we want show that
\begin{equation}\label{T10a}
\partial_t  \phi(\tov,\xo)+   \Ho_1(\xo,\overline  u(\tov, \xo),D\phi(\tov,\xo), D^2\phi(\tov,\xo))\le 0.
\end{equation}
Let  $v=v(y)$ be a periodic  viscosity solution of
\[
H_1(\xo,y,\overline  u(\tov,\xo), D\phi(\tov,\xo),D^2\phi(\tov,\xo) +D^2v(y))=\Ho_1(\xo,\overline  u(\tov,\xo),D\phi(\tov,\xo), D^2\phi(\tov,\xo));\]
    namely, $v$ solves the cell problem \eqref{CP} with $(x,r,p,X)=(\xo,\uo(\tov,\xo), D\phi(\tov,\xo), D^2\phi(\tov,\xo))$ (we recall that its existence is ensured by Proposition~\ref{Prp:erg}-(i)).
By \cite[Lemma 2.7]{hi} (recalled in Lemma \ref{lemma27hi} below) for each $\eta>0$, there exists a periodic supersolution $w\in C(\R^n)\cap W^{2,\infty}(\R^n)$ of
\begin{equation}\label{TH11}
    H_1( \xo,y,\overline  u(\tov,\xo), D\phi(\tov,\xo),D^2\phi(\tov,\xo) +D^2w(y))=\Ho_1(\xo,\overline  u(\tov,\xo),D\phi(\tov,\xo), D^2\phi(\tov,\xo))-\eta.
\end{equation}
Define the ``perturbed test-function''
\[
   \phi^\e(t,x)=\phi(t,x)+\e^2 w\left(\frac{x}{\e}\right).
\]
By standard results, we have that, up to extract subsequences, there exist
$(t_\e,x_\e)\in Q_T$, $(t_\e,x_\e)\to (\tov,\xo)$ for $\e\to 0$ such that $ (t_\e,x_\e)$ is a local maximum of ${u}^\e_1(t,x)-\phi^\e(t,x)$ and $ \mathop{\rm lim}_{\e\to 0}{u}^\e_1(t_\e,x_\e )= \overline u_1 (\tov,\xo)$.\\
Assume for the moment that $w\in C^2(\R^n)$ so that $\phi^\e$ is an admissible test function for $u^\e_1$ at $(t_\e,x_\e)$. Then
\begin{equation}\label{TH13}
\partial_t  \phi(t_\e,x_\e)+ H_1\left(x_\e,\frac{x_\e}{\e},{u}^\e(t_\e,x_\e),
D\phi(t_\e,x_\e)+\e Dw\left(\frac{x_\e}{\e}\right), D^2\phi(t_\e,x_\e)+D^2w\left(\frac{x_\e}{\e}\right)\right)
\leq 0.
\end{equation}
Set $\d_\e:=u_1^\e(t_\e,x_\e)-\phi^\e(t_\e,x_\e)$. By the definition of $\uo$, up to a subsequence,  for $j\neq 1$ $ u^\e_j(t_\e,x_\e)\to \bar r_j$ with $\bar r_j\le \uo_j(\tov,\xo)$ . By \eqref{TH13} and $(H5)$,
\begin{align*}
 &0\ge \partial_t \phi(t_\e,x_\e)+ H_1\left(x_\e,\frac{x_\e}{\e},{u}^\e(t_\e,x_\e),
D\phi(t_\e,x_\e)+\e Dw\left(\frac{x_\e}{\e}\right), D^2\phi(t_\e,x_\e)+D^2w\left(\frac{x_\e}{\e}\right)\right)\\
&=  \partial_t \phi(t_\e,x_\e)+ H_1\left(x_\e,\frac{x_\e}{\e},\left(\phi^\e(t_\e,x_\e)+\d_\e, u^\e_2(t_\e,x_\e),\dots,u^\e_m(t_\e,x_\e)\right),
D\phi(t_\e,x_\e)+\e Dw\left(\frac{x_\e}{\e}\right),\right.\\
&\qquad \left. D^2\phi(t_\e,x_\e)+D^2w\left(\frac{x_\e}{\e}\right)\right)\\
&\ge \partial_t  \phi(t_\e,x_\e)+ H_1\left(x_\e,\frac{x_\e}{\e},\left(\phi^\e(t_\e,x_\e),u^\e_2(t_\e,x_\e),\dots,u^\e_m(t_\e,x_\e)\right),
D\phi(t_\e,x_\e)+\e Dw\left(\frac{x_\e}{\e}\right),\right.\\
& \qquad\left. D^2\phi(t_\e,x_\e)+D^2w\left(\frac{x_\e}{\e}\right)\right)+\g \d_\e.
\end{align*}
We denote by $\xi$ the limit in $\R^n/\Z^n$ of $x_\e/\e$ as $\e \to 0$.
Passing to the limit for $\e\to 0$ in the previous inequality, by the  periodicity of $H_1$ and $w$,  \eqref{TH11}  and  $(H5)$
with $r=(\uo_1(\xo,\tov),\bar r_2,\dots,\bar r_m)$ and $s=(\uo_1(\tov,\xo),\uo_2(\tov,\xo),\dots,\uo_m(\tov,\xo))$ we get
\begin{eqnarray*}
0&\ge& \partial_t\phi(\tov,\xo)+ H_1(\xo,\xi, (\uo_1(\tov,\xo),\bar r_2,\dots,\bar r_m),
D\phi(\tov,\xo), D^2\phi(\tov,\xo)+D^2w(\xi))\\
&\ge&\partial_t \phi(\tov,\xo)+ H_1(\xo,\xi, (\uo_1(\tov,\xo),\uo_2(\tov,\xo),\dots,\uo_m (\tov,\xo)), D\phi(\tov,\xo), D^2\phi(\tov,\xo)+D^2w(\xi))\\
&\ge&\partial_t \phi(\tov,\xo)+\Ho_1(\xo,\uo(\tov,\xo),D\phi(\tov,\xo), D^2\phi(\tov,\xo))-\eta.
\end{eqnarray*}
and, for the arbitrariness of $\eta$ we get \eqref{T10a}.
If $w$ is not smooth, using  in \cite[Lemma 3.6]{hi} (recalled in Lemma \ref{lemma36hi}) it is possible to find $X_\e\in S^n$ such that
\begin{eqnarray*}
    (Dw(\frac{x_\e}{\e}), X_\e)&\in&\bar J^2w(\frac{x_\e}{\e})\\
    (D\phi(t_\e,x_\e)+\e Dw(\frac{x_\e}{\e}), D^2\phi(t_\e,x_\e)+X_\e)&\in& J^{2,+}u^\e(t_\e,x_\e)
\end{eqnarray*}
hence the above arguments 
 hold  with $X_\e$ in place of $D^2w(\frac{x_\e}{\e})$.
The rest of the proof to obtain \eqref{T10a} is exactly the same.\par
We prove that $\underline u$ is a viscosity supersolution of \eqref{hom} in a
similar way. From Proposition \ref{homo}, we then obtain   $\overline u\leq
\underline u$ in $Q_T$, hence  $\overline u=\underline u:=u$
where $u$ is the (local) uniform limit of the $u^\e$'s.
\end{Proof}
\begin{Remark}
Observe that in the previous proof we exploit three facts
\begin{itemize}
\item  for each $i\in I$, $H_i$ is ergodic, i.e. the cell problem \eqref{CP} admits a solution for any $(x,r,p,X)$.
\item there exist ``sufficiently  regular'' approximations to the solution to the cell problem \eqref{CP}
\item The effective Hamiltonian $\Ho_i$ satisfies the properties in Prop. \ref{Prp:erg}.ii).
\end{itemize}
The uniform ellipticity of $H_i$ is a sufficient condition to ensure these properties (for the last one, some regularity assumptions on the coefficients is also needed).
Let us stress that such properties still hold under different hypotheses as, for instance, for first order equations, the coercivity with respect to $p$ (in this case, the regular approximations of the solution to the cell problem will belong to $W^{1,\infty}$).
\end{Remark}
%
%
\begin{Remc}{Example 4.1}
Consider the weakly coupled system
\begin{equation}\label{lin}
\partial_t u^\e_i-\textrm{tr} \left(a_i\big(x,\frac x\e\big) D^2u^\e_i\right)+F_i\left(x,\frac x\e,u,Du^\e_i\right)=0\qquad (t,x)\in Q_T,\,i\in I
\end{equation}
where
$F_i\left(x,y,r,p\right)=\min\limits_{\z\in Z}\max\limits_{\th\in \Theta}\left\{
-f_i^{\tz}(x,y)\cdot p-l_i^{\tz}(x,y) -\sum_{j}d^{\tz}_{ij}(x,y)r_j\right\}$.

For each $\xo,\uo,\po,\Xo$, the cell problem reads
\begin{equation*}
-\textrm{tr} \left(a_i(\xo,y) D^2_y v\right)-\textrm{tr} \left(a_i(\xo,y)\Xo\right)+F_i\left(\xo,y,\uo,\po\right)=\Ho(\xo,\uo,\po,\Xo).
\end{equation*}
By standard theory for linear ergodic problems (see \cite{BLP} and also \cite{AB8}), there holds
\begin{equation}\label{lineff}
\Ho(x,u,p,X)=-\textrm{tr} \left(\bar a_i(x) X\right)+\overline F_i\left(x,u,p\right)
\end{equation}
where the effective diffusion $\bar a$ and the effective operator $\bar F_i$ have respectively the form
$$
\bar a_i(x):=\int_{[0,1)^n} a_i(x,y)\, d\mu_x(y),\quad
\overline F_i\left(x,r,p\right):=\int_{[0,1)^n} F_i\left(x,y,r,p\right)\, d\mu_x(y).
$$
Here, for $\xo$ fixed, the measure $\mu_{\xo}$ is the unique invariant measure for the diffusion $a(\xo,y)$,
i.e. the solution in the sense of distributions of the equation
$$
\sum_{i,j=1}^n \frac{\partial^2}{\partial y_i \,\partial y_j}\left(a_{ij}(\xo,y)\mu_{\xo}\right)=0,\qquad \mu_{\xo} \textrm{ periodic.}
$$

As a straightforward application of Proposition \ref{Prp:erg}-(ii), Proposition \ref{homo} and Theorem \ref{thm:hom}, we have the following result
\begin{Corollary}
Let $u^\e$ and $u$ be respectively the solution to system \eqref{HJe} with $H_i$ as in \eqref{lin} and the solution to \eqref{hom} with $\Ho$ as in \eqref{lineff}. Then $u^\e$ converges locally uniformly to $u$ on $[0,T]\times \R^n$.
\end{Corollary}
\end{Remc}

%
%
\appendix
\section{Appendix}\label{appendix}
For the proof of Lemma \ref{Lemma2}, we need  the following technical Lemma:
\begin{Lemma}\label{Lemma1}
Let $f\in USC(\R^N\times \R_+\times I)$ be bounded from above and $g\in C(\R^N\times \R_+)$ be nonnegative. For $\e>0$, set $\psi_\e(\xi,t,i):=f(\xi,t,i)-\e g(\xi,t)$ and assume that $\psi_\e$ attains its global maximum in some point $(\xi_0^\e,t_0^\e,i_0^\e)$. Then, as $\e \to 0$, $\max \psi_\e \to \sup f$ and $\e g(\xi_0^\e,t_0^\e)\to 0$.
\end{Lemma}
\begin{Proof} Set $m_\e:=\max \psi_\e$ and $m:=\sup f$. For $\eta>0$, let $(\xi',t',i')$ be such that: $f(\xi',t',i')\ge m-\eta$. For $\e'$ sufficiently small, we have: $\e' g(\xi',t')\le \eta$. In particular, since $g$ is nonnegative, there holds
\begin{equation*}
m\ge m_{\e'}\ge f(\xi',t',i')-\e' g(\xi',t')\ge m-2\eta.
\end{equation*}
Letting $\e'\to 0$, we get the first part of the statement.

For $\e$ sufficiently small, the above relations entail
$$
m_\e= f(\xi_0^\e,t_0^\e,i_0^\e)-\e g (\xi_0^\e,t_0^\e)\ge m-2\eta;
$$
in particular, for $k_\e:=\e g (\xi_0^\e,t_0^\e)$, we deduce that the sequence $\{k_\e\}_\e$ is bounded.
Let us pick a subsequence (still denoted $k_\e$) convergent to some value $k\ge0$.
Since $m_\e=f(\xi_0^\e,t_0^\e,i_0^\e)-k_\e\le m-k_\e$, by the first part of the statement, as $\e \to0$, we obtain $k\le0$. Hence $k=0$ and the statement is completely proved.
\end{Proof}

\begin{Proofc}{Proof of Lemma~\ref{Lemma2}}
$(i)$. Relations~\eqref{T4} and~\eqref{T5} entail
\begin{align*}
    0\le \psi(\t_0,x_0,y_0,i_0)\le 2R
-\left(\frac \a 2 e^{\bar \g \t_0}|x_0-y_0|^2+\frac \e 2(|x_0|^2+|y_0|^2)\right);
\end{align*}
therefore, inequalities \eqref{T7} easily follows. The estimates~\eqref{T8} are an immediate consequence of \eqref{T7} and Lemma~\ref{Lemma1}.

$(ii)$. The inequality $2\psi(\t_0,x_0,y_0,i_0)\geq \psi(\t_0,x_0,x_0,i_0)+\psi(\t_0,y_0,y_0,i_0)$ yields
$$
\a e^{\bar \g \t_0}|x_0-y_0|^2\leq \left[u^1_{i_0}(\t_0,x_0)-u^1_{i_0}(\t_0,y_0)\right]+\left[u^2_{i_0}(\t_0,x_0)-u^2_{i_0}(\t_0,y_0)\right];
$$
therefore, inequality \eqref{T7-1} is a consequence of the regularity assumption.

$(iii)$. Assume $u^1\in C^1$ (being similar, the other case will be omitted). Let $\{e_k\}$ be an orthogonal basis of $\Rn$. For $h\in\R$ sufficiently small, the inequality $\psi(\t_0,x_0,y_0,i_0)\geq \psi(\t_0,x_0+he_i,y_0,i_0)$ yields
$$
\a e^{\bar \g \t_0}\left(|x_0-y_0+h e_k|^2-|x_0-y_0|^2\right) +\frac \e2\left(|x_0+h e_k|^2-|x_0|^2\right)\leq u^1_{i_0}(\t_0,x_0)-u^1_{i_0}(\t_0,x_0+he_k).
$$
Dividing by $h$ and letting $h\to 0^\pm$, we obtain
$$
|\a e^{\bar \g \t_0}(x_{0,k}-y_{0,k})+\e x_{0,k}|\leq [u^1]_1.
$$
Summing on $k$ and taking advantage of estimate~\eqref{T7}, we conclude the proof.
\end{Proofc}

For the sake of completeness, let us now state two results established by Horie and Ishii in~\cite[Lemma 2.7 and 3.6]{hi}. For their proof, we refer the reader to the original paper.

\begin{Lemma}\label{lemma27hi}
Assume conditions $(H0)$-$(H4)$ and fix $x,p\in\R^n, r\in\R^m, X\in S^n, i\in I$. Let $v=v(y)$ be a bounded continuous solution to \eqref{CP}. Then
\begin{itemize}
\item[(a)] $v$ is Lipschitz continuous in $\R^n$.
\item[(b)] Let $R>0$ be a constant such that $\|Dv\|\leq R$. Then, for each $\e>0$, there are functions $v^\pm\in C(\R^n)\cap W^{2,\infty}(\R^n)$ and a constant $C$ (depending on $R$ and on the constants entering in the assumptions) such that
\begin{align*}
&\|v-v^\pm\|\leq \e,&\qquad \|v^\pm\|\leq \|v\|\\
&\|D v^\pm\|\leq \|Dv\|, &\qquad \|v^\pm\|_{W^{1,\infty}(\R^n)}\leq C,
\end{align*}
and
\begin{eqnarray*}
 H_i(x,y,r,p, X+D^2v^+(y))&\geq&\Ho_i(x,r,p,X)-\e \qquad \textrm{in }\R^n\\
 H_i(x,y,r,p, X+D^2v^-(y))&\leq&\Ho_i(x,r,p,X)+\e \qquad \textrm{in }\R^n.
\end{eqnarray*}
\end{itemize}
\end{Lemma}

\begin{Lemma}\label{lemma36hi}
Let $\Omega\subset \R^n$ be open, $u\in USC(\Omega)$ and $v\in C(\Omega)\cap W^{2,\infty}(\Omega)$. Let $\hx\in \Omega$ and $(p,X)\in J^{2,+}(u-v)(\hx)$. Then there exists a $Y\in S^n$ such that
\begin{equation*}
(Dv(\hx),Y)\in J^2 v(\hx),\qquad
(p+Dv(\hx), X+Y)\in J^{2,+}u(\hx)
\end{equation*}
where $J^{2,+}u(\hx)$ is the set of superjets of $u$ at the point $\hx$ (see \cite[Section 2]{CIL}) while $\bar J ^2v(x)$ denotes the set of those points $(q,Y)\in \R^n\times S^n$ for which there is a sequence $x_j\to x$ such that $v$ is twice differentiable at $x_j$ and $(Dv(x_j), D^2v(x_j))\to (q,Y)$ (see \cite[Section 3]{CIL}).
\end{Lemma}

\begin{Remc}{\bf Acknowledgement.} The authors are grateful to Professor E.R. Jakobsen for several useful comments and suggestions.

\end{Remc}

\end{document}